\documentclass[dvips,preprint,authoryear,12pt]{imsart}

\RequirePackage{natbib,amssymb,verbatim,undertilde,mathrsfs}
\RequirePackage[OT1]{fontenc}
\RequirePackage{amsthm,amsmath,natbib}
\RequirePackage{hypernat}
\usepackage{graphicx}
\def\newblock{\hskip .11em plus .33em minus .07em}

\def\a{\alpha}

\def\bv{\mathbf{v}}

\def\1{\mathbf{1}}

\def\argmin{\operatorname{argmin} \displaylimits}
\def\l{\lambda}

\def\e{\epsilon}
\def\ee{\boldsymbol{\epsilon}}
\def\b{\beta}

\def\bb{\boldsymbol{\beta}}

\def\mm{\boldsymbol{\mu}}

\def\x{\mathbf{x}}
\def\w{\mathbf{w}}

\def\s{\sigma}

\def\u{\mathbf{u}}
\def\y{\mathbf{y}}

\def\1{\mathbf{1}}

\def\Dto{\stackrel{D}{\to}}

\def\R{\mathbb{R}}

\def\W{\mathbf{W}}
\def\ss{\mathbf{s}}
\def\rr{\mathbf{r}}

\def\X{\mathbf{X}}

\begin{document}

\bibliographystyle{ims}
\begin{frontmatter}
\title{Parallelism, uniqueness, and large-sample asymptotics for the Dantzig selector \protect}
\runtitle{Uniqueness and asymptotics for the Dantzig selector}

\begin{aug}
\author{\fnms{Lee} \snm{Dicker}%\thanksref{t1}
\ead[label=e1]{ldicker@stat.rutgers.edu}}
\and \author{\fnms{Xihong} \snm{Lin}%\thanksref{t2}
\ead[label=e2]{xlin@hsph.harvard.edu}}

%\thankstext{t1}{Supported by NIH Grant T32 ES007142}

\runauthor{L. Dicker and X. Lin}

\affiliation{Rutgers University}

\address{Department of Statistics and Biostatistics \\ Rutgers University \\ 501 Hill Center, 
 110 Frelinghuysen Road \\ Piscataway, NJ 08854 \\
\printead{e1}}
\address{Department of Biostatistics \\ Harvard School of Public
  Health \\ 655 Huntington Avenue \\ Boston, MA 02115 \\
\printead{e2}}
\end{aug}


\begin{abstract}
The Dantzig selector \citep{candes2007dantzig} is a popular
$\ell^1$-regularization method for variable selection and estimation
in linear regression.  We present a very weak geometric condition on
the observed predictors which is related to parallelism and, when satisfied,
ensures the uniqueness of Dantzig selector
estimators. The condition holds with
probability 1, if the predictors are drawn from a continuous
distribution.  We discuss the necessity of this condition for
uniqueness and also provide a closely
related condition which ensures uniqueness of lasso estimators
\citep{tibshirani1996regression}. Large sample
asymptotics for the Dantzig selector, i.e. almost sure convergence and the asymptotic
distribution, follow directly from our uniqueness results and a
continuity argument.  The limiting
distribution of the Dantzig selector is generally non-normal. Though our
asymptotic results require that the number of predictors is fixed
(similar to \citep{knight2000asymptotics}), our uniqueness results are
valid for an arbitrary number of predictors and observations.
\end{abstract}

\end{frontmatter}

\section{Introduction}

Regularized regression methods for variable selection and estimation
have become an important tool for statisticians
and have been the subject of intense statistical research during the past
fifteen years \citep{bickel2006regularization,fan2010selective,tibshirani2011regression}.  These methods
provide a tractable approach to the
analysis of high-dimensional datasets and are especially useful when
the underlying signal is sparse.  In this paper, we address some gaps
in the literature, which pertain to uniqueness and large sample
asymptotic theory for the Dantzig selector \citep{candes2007dantzig},
a popular $\ell^1$-regularized regression method that is closely related to lasso
\citep{tibshirani1996regression}.

First, we develop an intuitive geometric condition related to parallelism which ensures that
the Dantzig selector has a unique solution and demonstrate that
this condition holds in an overwhelming majority of instances (with
probability 1, if the predictors follow an absolutely continuous
distribution with respect to Lebesgue measure).  We also give a
related necessary condition for the uniqueness of Dantzig selector
solutions.  These results originally appeared in the first
author's PhD thesis \citep{dicker2010thesis} and, to our knowledge, are the first
uniqueness results about the Dantzig selector to be found in the
literature.  In fact, our uniqueness condition for the Dantzig
selector is easily translated into a similar prevalent condition which
implies that lasso has a unique solution.

 Aside from their independent interest, the uniqueness results
 presented here pave the way for a simple derivation of the almost sure limit and the asymptotic distribution of Dantzig
selector estimators, when the number of predictors, $p$, is fixed (on
the other hand, we
emphasize that our uniqueness results are valid for arbitrary $p$).  These
asymptotic results are analogous to those found in
\citep{knight2000asymptotics}  for the lasso and further highlight
similarities between the two methods, which have been discussed by
multiple authors \citep{meinshausen2007discussion, james2009dasso}.  
In fact, in comparison with Knight and Fu's [2000] results, uniqueness
appears to be the major hurdle to obtaining large sample asymptotics
for the Dantzig selector.  The Dantzig selector is a convex -- but not strictly convex
-- optimization problem.  Thus, unique solutions are not guaranteed in
general.  However, once uniqueness is understood, asymptotic results
for the Dantzig selector follow directly from continuity arguments.
More specifically, we show that under the given uniqueness
conditions the Dantzig selector may be viewed as a well-defined
continuous mapping; asymptotic results then follow from the
continuous mapping theorem.   By contrast, for the lasso, uniqueness
is assured in classical fixed $p$ asymptotic analyses because the
associated optimization problem is strictly convex (provided the
predictors are non-degenerate).  The foregoing discussion highlights the potential
usefulness of uniqueness results for the Dantzig selector.  More
broadly, understanding uniqueness makes certain powerful tools -- like the continuous mapping theorem -- readily available for further analysis
of the Dantzig selector.  

Though much of the recent
interest in regularized regression methods is spurred by applications
that may perhaps be best approximated by an asymptotic regime where
$p \to \infty$, we believe that it remains important to understand
classical large sample asymptotics, where $p$ is fixed
  and $n \to \infty$, in order to obtain a more complete understanding
  of these procedures.  This paper helps shed light on this issue.   Moreover, we believe
that our uniqueness results, which are valid for all $p$, may be
useful for formulating and deriving asymptotic results for regularized
regression methods in settings where $p \to \infty$; however, this is
a topic for future research and is beyond the scope of this paper (though
it is briefly addressed again in our concluding Section 5).  

The rest of this paper proceeds as follows.  In Section 2 we introduce
notation and definitions.  In Section 3 we discuss uniqueness.  Propositions 1 and 2 are
the main results in Section 3 and summarize important uniqueness properties of
the Dantzig selector and lasso vis-\`a-vis parallelism.  In
Section 4, we show that the Dantzig selector may be viewed as a
continuous mapping from the space of predictors and associated
outcomes to the space of parameter estimates (Proposition 3).
 Corollaries 1 and 2 give the almost-sure limit of
 Dantzig selector estimators and their asymptotic distribution,
 respectively.  Section 5 contains a brief concluding discussion.  Proofs may be found in the
Appendix at the end of the paper.

\section{Notation and definitions}

Consider the linear model
\begin{equation} \label{lm1}
y_i = \x_i^T\bb^* + \e_i, \ \ i = 1,...,n
\end{equation}
where $y_1,...,y_n \in \R$ and $\x_1,...,\x_n \in \R^p$ are observed
outcomes and predictors, respectively, $\e_1,...,\e_n$ are unobserved
iid integrable random variables with mean $E(\e_i) = 0$, and $\bb^* = (\b_1^*,...,\b_p^*)^T \in \R^p$ is an
unknown parameter to be estimated.  To simplify notation, let $\y = (y_1,...,y_n)^T \in \R^n$
denote the $n$-dimensional vector of outcomes and $X =
(\x_1,...,\x_n)^T$ denote the $n\times p$ matrix of predictors.  Also
let $\ee = (\e_1,...,\e_n)^T \in \R^n$.   Then (\ref{lm1}) may be
re-expressed as
\[
\y = X\bb^* + \ee.
\]

It will be useful to have a concise method for referring to
sub-vectors and sub-matrices of various vectors and matrices.  For a
vector $\bb = (\b_1,...,\b_p)^T\in \R^p$ and a subset $A \subseteq \{1,...,p\}$, let $\bb_A
= (\b_j)_{j \in A} \in \R^{|A|}$.  Furthermore, for $n \times p$
matrices $X = (x_{ij})_{1 \leq i \leq n, \ 1 \leq j \leq p}$ let $X_A =
(x_{ij})_{1 \leq i \leq n, \ j \in A}$ denote the $n \times |A|$ matrix
obtained from $X$ by extracting columns corresponding to elements of
$A$.  If $C = (c_{ij})_{1 \leq i,j \leq p}$ is a $p \times p$ matrix, and $B \subseteq \{1,...,p\}$
has cardinality $|B|$, let $C_{A,B} = (c_{ij})_{i  \in A, \ j \in B}$ denote the $|A|
\times |B|$ matrix obtained from $C$ by extracting rows corresponding
to elements of $A$ and columns corresponding to elements of $B$.
For $j \in \{1,...,p\}$, let $\X_j =
X_{\{j\}}$ denote the $j$-th column of $X$.  Finally, let
$\mathrm{null}(C)$ denote the null-space of the matrix $C$ and let
$\dim(V)$ denote the dimension of the vector space $V$.

The main object of study in this paper is the Dantzig selector -- a
linear programming problem for obtaining estimates of $\bb^*$, which is
defined as follows:
\begin{equation}\label{DS}
\begin{array}{lc}
\mbox{minimize} & ||\bb||_1 \\
\mbox{subject to} & \frac{1}{n} ||X^T(\y - X\bb)||_{\infty} \leq \lambda,
\end{array}
\end{equation}
where $\lambda = \l_n \geq 0$ is a tuning parameter, $||\bb||_1 = \sum_{j =
  1}^p |\b_j|$ denotes the $\ell^1$-norm and $||X^T(\y -
X\bb)||_{\infty} = \max_{1 \leq j \leq p} |\X_j^T(\y - X\bb)|$ denotes the
$\ell^{\infty}$-norm.  Solutions to (\ref{DS}), denoted
$\hat{\bb}^{ds}$, will be referred to as
Dantzig selector estimators. 

We also introduce the lasso optimization problem and estimator
at this time:
\begin{equation}\label{lasso}
\hat{\bb}^{lasso} \in \argmin_{\bb \in \R^p} \frac{1}{2n} ||\y - X\bb||^2 +
\l ||\bb||_1,
\end{equation}
where $||\y - X\bb||^2 = \sum_{i = 1}^n (y_i - \x_i^T\bb)^2$ is the squared
$\ell^2$-norm.  Though the lasso is not our primary concern in this
paper, we will sometimes find it instructive to compare aspects of the
Dantzig selector and lasso side-by-side.  For instance, as discussed in the
Introduction, notice that if $X$ has rank $p$, then lasso is a
strictly convex optimization problem, which ensures that
$\hat{\bb}^{lasso}$ is unique.  On the other hand, the Dantzig selector
(\ref{DS}) is a linear programming problem and uniqueness properties
are less clear, even when $X$ has rank $p$.

In order to provide some additional context for the present study, we point out that one of the key
features of both the Dantzig selector
and lasso is that they perform simultaneous variable selection and
estimation.  By this we mean that $\{j; \ \hat{\b}^{ds}_j = 0\}$
and $\{j; \ \hat{\b}^{lasso}_j = 0\}$ are often non-empty
(contrast this with the ordinary least squares estimator for $\b^*$).
This implies that $\hat{\bb}^{ds}$ and
$\hat{\bb}^{lasso}$ often have reduced dimension (i.e., only a few
non-zero entries) and can greatly enhance interpretability, along with
estimation accuracy \citep{tibshirani1996regression, candes2007dantzig,
  bickel2009simultaneous}.

\section{Parallelism and uniqueness}

Parallelism plays a large role in the discussion of uniqueness of Dantzig selector solutions.  Roughly
speaking, the Dantzig selector has a unique solution if the feasible set,
\[
F = \{\bb; \ ||X^T(\y - X\bb)||_{\infty} \leq \l\} \subseteq \R^p,
\]
is not parallel to the $\ell^1$-ball.  Below, we describe parallelism as a geometric concept which is
relevant to the Dantzig selector and then give a more formal definition.

First note that the feasible set $F$ is polyhedral (it is the
intersection of finitely many hyperplanes).  Solutions of the
Dantzig selector are points $\b \in F$ of minimal $\ell^1$-norm.  Let
$B_1 = \{\u \in \R^p; \ ||\u||_1 \leq 1\}$ be the closed unit
$\ell^1$-ball centered at the origin.  Geometrically, we can find
solutions to the Dantzig selector by ``growing" $tB_1 = \{\u \in \R^p;
\ ||\u||_1 \leq t\}$, $t \geq 0$, until it intersects $F$; the points
of intersection are Dantzig selector solutions.  More precisely, let
$t_0 = ||\hat{\bb}^{ds}||_1$.  The collection of all Dantzig selector solutions is $F \cap
t_0B_1$.  When $p = 2$, the 1-dimensional faces of $tB_1$ have slope 1
or $-1$; the Dantzig selector has multiple solutions only if a
1-dimensional face of $F$ has slope 1 or -1, that is, only if $F$ is
parallel to the $\ell^1$-ball, $B_1$.

\begin{figure}[tbp]
\begin{center}
\hspace{-.3in}\includegraphics[width=0.7\textwidth]{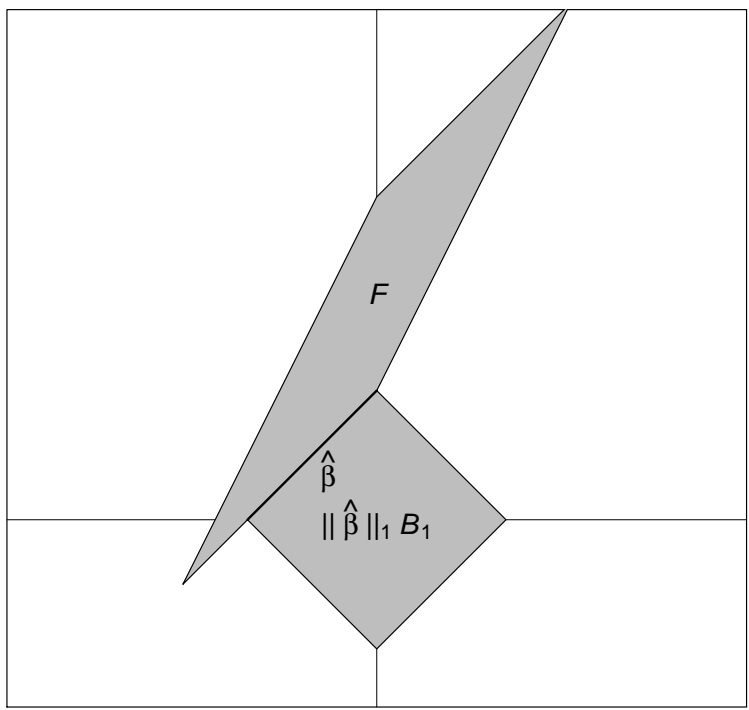}
\end{center}

%\vspace{-.5in}
\caption{An instance of the Dantzig selector with multiple solutions.  The region $F$ is the feasible set
for the Dantzig selector and $||\hat{\b}||_1B_1 = \{u\in \R^2; \ ||u||_1 \leq ||\hat{\b}||_1\}$.  The bold
line represents the intersection of $||\hat{\b}||_1B_1$ with $F$ and is the solution set for this instance
of the Dantzig selector.}
\end{figure}

As indicated by the situation when $p = 2$, if the Dantzig selector
has  multiple solutions, then $F$ is parallel to $B_1$ (Figure 1).  When $p \geq 2$, the  notion of
parallelism which is correct for our purposes is less straightforward.  Geometric intuition suggests that
parallelism is invariant under translation and scalar multiplication,
in the sense that $F$ is parallel to $B_1$ if and only if $\a F + \mathbf{v}_0
= \{\a \bb + \mathbf{v}_0; \ \bb \in F\}$ is
parallel to $B_1$ for $a \in \R\setminus\{0\}$ and $\mathbf{v}_0 \in
\R^p$.  In particular, multiplying $X$ by a (non-zero) scalar and
adding vectors $\y_0 \in \R^n$ to $\y$ does not affect parallelism.  This leads to a
definition of parallelism between $F$ and $B_1$ which depends only on the matrix $n^{-1}X^TX$.  In fact,
in our view, the primitive concept is parallelism between a $p\times p$ symmetric matrix $C$ and the
$\ell^1$-ball.

\vspace{.1in}
\noindent {\bf Definition 1. } {\em \begin{itemize} \item[{\em{\bf (a)}}] Let $C$ be a $p \times p$ symmetric matrix.  The matrix $C$ is
parallel to the $\ell^1$-ball if and only if the condition [Par] (found below) holds.
\begin{list}{[Par]}
\item There exist subsets $A,B \subseteq \{1,...,p\}$ and a vector $\w \in
  \R^{|B|}$ such that $||C_B \w||_{\infty} \leq 1$, $C_{A,B}\w \in
  \{\pm 1\}^{|A|}$, and
  $\dim\left[\mathrm{null}\left(C_{B,A}\right)\right]  > 0$.
\end{list}

\item[{\em{\bf (b)}}]  The feasible set for the Dantzig selector, $F$, is parallel to the $\ell^1$-ball if and
only if $n^{-1}X^TX$ is parallel to the $\ell^1$-ball. \end{itemize}}

\vspace{.1in}

\noindent {\em Remarks (i)} Parallelism, as defined here, is related
to degenerate sub-matrices of $C$, which, in the context of the
Dantzig selector, correspond to the nontrivial faces of $F$.  In [Par], the requirement
that $C_{A,B}w \in \{\pm 1\}^{|A|}$ is related to the fact that the
faces of the $\ell^1$-ball, $B_1$, have normal vectors $u \in \R^p$,
where $u_A \in \{\pm1\}^{|A|}$ for some $A \subseteq \{1,...,p\}$.

\noindent {\em (ii)} When $p = 2$, it is easy to see that $F$ is parallel to the $\ell^1$-ball if and only
if one of the columns of $n^{-1}X^TX$ is a scalar multiple of some point in $\{\pm1\}^2$.  This occurs if
and only if a one-dimensional face of $F$ has slope 1 or -1, as depicted in Figure 1.

%\vspace{-.2in}

\begin{figure}[tbp]
\begin{center}
\includegraphics[width = \textwidth]{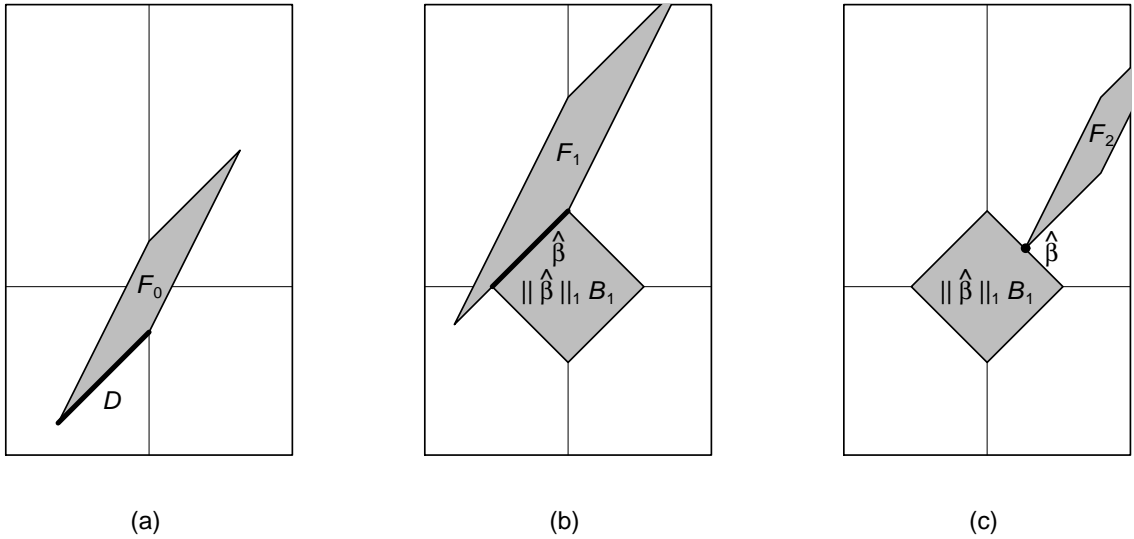}
\end{center}
\vspace{-.2in}
\caption{(a) $F_0 = (X'X)^{-1}B_{\infty}$ is parallel to the $\ell^1$-ball, as evidenced by the bold face
$D$.  (b) $F_1$ is obtained from $F_0$ by scalar multiplication and translation; the Dantzig selector
problem with feasible set $F_1$ has multiple solutions, indicated by the bold line segment labeled
$\hat{\b}$.  (c) $F_2$ is obtained from $F_0$ by scalar multiplication and translation; the point labeled
$\hat{\b}$ is the unique solution to the Dantzig selector problem with feasible set $F_2$.}
\end{figure}

As discussed above, parallelism is invariant under translation and
scalar multiplication.  On the other hand, translation and scalar
multiplication of the feasible set $F$ gives rise to various instances
of the Dantzig selector, some with a unique solution and some,
perhaps, with multiple solutions.  This suggests that any sufficient
condition for the existence of multiple Dantzig selector solutions
must, unlike parallelism, involve $\y$ and $\l$.  To illustrate this
concept, suppose that $n^{-1}X^TX$ is invertible and is parallel to
the $\ell^1$-ball.  Figure 2 (a) depicts $F_0 =
(n^{-1}X^TX)^{-1}B_{\infty} = \{(n^{-1}X^TX)^{-1}\u; \ \u \in \R^2,\
||\u||_{\infty} \leq 1 \}$, which is equal to the feasible set for the
Dantzig selector when $\l = 1$ and $\y = 0$ and is parallel to the
$\ell^1$-ball.  Figures 2 (b) and (c) depict $F_1$ and $F_2$,
potential feasible sets for the Dantzig selector that are both
obtained from $F_0$ by scalar multiplication and translation.  The
feasible sets $F_1$ and $F_2$ are both parallel to the $\ell^1$-ball,
and correspond to feasible sets for the Dantzig selector with the
predictor matrix $X$ and different values for $\y$, $\l$ (not given
here).  The instance of the Dantzig selector with feasible set $F_1$
has multiple solutions, while the Dantzig selector with feasible set
$F_2$ has a unique solution.

The following condition combines parallelism with additional constraints and is a sufficient condition for
the existence of multiple Dantzig selector solutions.
\begin{description}
\item[{\normalfont [Mult]}]{\em There exist subsets $A,B \subseteq \{1,...,p\}$ and vectors $\mm^0 \in \R^{|B|}$, $\bb^0 \in F$, such that
\begin{itemize}
\item[{\em 1.}] $||n^{-1}X^TX_B\mm^0||_{\infty} \leq 1$, $n^{-1}X_A^TX_B\mm^0 \in \{\pm
  1\}^{|A|}$, and
  $\dim\left[\mathrm{null}\left(n^{-1}X_B^TX_A\right)\right]  > 0$.
\item[{\em 2.}] $n^{-1}\bb^TX^TX_B\mm^0 \geq ||\bb^0||_1$ for all $\bb \in F$.
\item[{\em 3.}] $A = \{j; \ \b_j^0 \neq 0\}$.
\item[{\em 4.}] $n^{-1}|\X_j^T(\y - X\bb^0)| = \l$ for all $j \in B$ and $n^{-1}|\X_j^T(\y - X\bb^0)| < \l$ for
    all $j \notin B$.
\end{itemize}}
\end{description}

Note that Condition 1 in [Mult] implies that $F$ is parallel to the
$\ell^1$-ball.  Conditions 2-4 in [Mult] constrain the location of $F$
in $\R^p$ relative to the origin.  Proposition 1 below characterizes
uniqueness properties of the Dantzig selector in terms of [Par] and
[Mult].  A related necessary condition for the existence of multiple lasso
solutions is given in Proposition 1 (c).  Proposition 1
is proved in the Appendix at the end of this paper.

\vspace{.1in}
\noindent {\bf Proposition 1. }{\em\begin{itemize} \item[{\em{\bf (a)}}] If [Mult] holds, then the Dantzig selector has multiple solutions.
\item[{\em {\bf (b)}}] If $F$ is not parallel to the $\ell^1$-ball, then the Dantzig
selector has a unique solution.
\item[{\em {\bf (c)}}] Suppose that $\l > 0$ and that the lasso has multiple solutions
(i.e. $\argmin_{\bb \in \R^p} (2n)^{-1}||\y - X\bb||^2 + \l||\bb||_1$
contains more than a single element).  Then there exists a subset $A \subseteq \{1,...,p\}$
and a vector $w \in \R^p$ such that $||n^{-1}X^TX\w||_{\infty} \leq 1$, $n^{-1}X_A^TX\w
\in
\{\pm1\}^{|A|}$, and
$\dim\left[\mathrm{null}\left(n^{-1}X^TX_A\right)\right] > 0$. \end{itemize}}

\vspace{.1in}

\noindent {\em Remarks (i) }  Proposition 1 is valid for any $n$ and
$p$.

\noindent {\em (ii) } Proposition 1 (c) may be rephrased as follows.
If the lasso has multiple solutions, then $n^{-1}X^TX$ is parallel to
the $\ell^1$-ball and, moreover, one may take $B = \{1,...,p\}$ in the
definition of parallelism.

\noindent {\em (iii) } If $\l =  0$, then the lasso has multiple
solutions whenever $n^{-1}X^TX$ is singular.

\noindent {\em (iv) } The condition in Proposition 1 (c) implies that $n^{-1}X^TX$ is
parallel to the $\ell^1$-ball.  It follows that if
$n^{-1}X^TX$ is {\em not} parallel to the $\ell^1$-ball, then both the
Dantzig selector and lasso have unique solutions.
The relationship between uniqueness for the Dantzig selector and
uniqueness for lasso is discussed by \cite{meinshausen2007discussion},
who give a concrete $p=3$-dimensional example (with pictures) where lasso has a unique solution, but the
Dantzig selector does not.  

\noindent {\em (v) } A condition similar to [Mult] which ensures the
existence of multiple lasso solutions may be developed.  This is not
pursued further here .

The next proposition suggests that the Dantzig selector and lasso
have a unique solution in an overwhelming majority of instances.

\vspace{.1in}
\noindent {\bf Proposition 2. }{\em
Suppose that $\x_1,...,\x_n$ are iid and drawn from a continuous
distribution with respect to Lebesgue measure on $\R^p$.  Then
$n^{-1}X^TX$ is parallel to the $\ell^1$-ball with probability 0.  Consequently,
the Dantzig selector and lasso have a unique solution with probability 1.}

\vspace{.1in}

\noindent {\em Remarks (i) } Proposition 2 is proved in
the Appendix (a proof also appears in \citep{dicker2010thesis}).
To provide some intuition, note that the
parallelism condition requires $n^{-1}X^T_AX_B$ to both (i) contain a
specific point in its range (that is, an element of $\{\pm 1\}^{|A|}$) and
(ii) to have a degenerate range (in the sense that
$\dim\left[\mathrm{null}(n^{-1}X^T_BX_A)\right] > 0$).  
Proposition 2 implies that this occurs with probability 0, under
the specified conditions.

\section{Large sample asymptotics for the Dantzig selector}

Throughout the rest of this article, assume that $p$ and $\bb^* \in
\R^p$ are fixed.  In this section, we formulate the Dantzig selector as a well-defined
mapping from sample covariance matrices, $n^{-1}X^TX$,
 marginal covariances, $n^{-1}X^T\y$, and tuning parameters, $\l \geq
 0$, to estimators, $\hat{\bb}^{ds}$.  To do this, we restrict our attention
to symmetric matrices that are {\em not} parallel to the
$\ell^1$-ball -- Proposition 2 suggests that this restriction is fairly
weak. Then, we show that the Dantzig selector mapping is continuous.
With this machinery in place, large
sample asymptotics for the Dantzig selector follow easily.

Let $\mathscr{P}_0$ denote the collection of $p \times p$ positive semidefinite
matrices that are not parallel to the $\ell^1$-ball and let
$\mathscr{P}_0^+ = \mathscr{P}_0 \cap \mathrm{GL}(p)$, where
$\mathrm{GL}(p)$ is the collection of all invertible $p \times p$
matrices with real entries.  Define the Dantzig selector mapping $G: \mathscr{P}_0^+ \times \R^p \times \R^{\geq 0} \to \R^p$
by $G(C,\bv,\lambda) = \hat{\u}$, where $\hat{\u}$ solves the optimization problem
\begin{equation}\label{defG}
\begin{array}{lc}
\mbox{minimize} & ||\u||_1 \\
\mbox{subject to} & ||C\u - \bv||_{\infty} \leq \lambda.
\end{array}
\end{equation}
It follows directly from Proposition 1 (b) that $G$ is well-defined.
Furthermore, notice that $G(n^{-1}X^TX,n^{-1}X^T\y,\lambda) =
\hat{\bb}^{ds}$.  Note that the domain of $G$ may be extended to a subset of
$\mathscr{P}_0 \times \R^p \times \R^{\geq 0}$, provided one imposes
conditions to ensure that the feasible set in the optimization problem (\ref{defG})
is non-empty.  More specifically, define $\mathscr{Q} = \{(C,\mathbf{v}); \ C
\in \mathscr{P}_0, \ \mathbf{v} \in \mathrm{range}(C)\}$.  Then (\ref{defG}) defines
$G(C,\mathbf{v},\l)$ for $(C,\mathbf{v},\l) \in \mathscr{Q} \times \R^{\geq 0}$.  

\vspace{.1in} 

\noindent {\bf Proposition 3. } {\em The mapping $G$ is continuous on $\mathscr{P}_0^+\times \R^p \times \R^{\geq 0}$.}

\vspace{.1in}

\noindent {\em Remarks (i) } A proof of
Proposition 3 is found in the Appendix.  A similar proof shows that $G$ is also
continuous on $\mathscr{Q} \times \R^{>0}$.  In other
words, assuming that the appropriate (anti-) parallelism conditions
hold, if there is non-trivial regularization in the limit (i.e. $\l_n
\to \l_0 > 0$), then the Dantzig selector is continuous, regardless of whether or not the
predictors and the limiting sample covariance matrix are singular.

\vspace{.1in}
\noindent {\bf Corollary 1. } {\em
Suppose that $n^{-1}X^TX \to C \in \mathscr{P}_0^+$ and that $\l_n \to \l_0 \geq 0$.  Then $\hat{\bb}^{ds}
\to \bb^0$, almost surely, where $\bb^0$ solves
\[
\begin{array}{lc}
\mbox{minimize} & ||\bb||_1 \\
\mbox{subject to} & ||C(\bb - \bb^*)||_{\infty} \leq \l_0.
\end{array}
\]}

\vspace{.1in}

\noindent {\em Remarks (i)} The corollary follows directly from
Proposition 3, which implies that $\hat{\bb}^{ds} = G(n^{-1}X^TX,n^{-1}X^T\y,\l_n) \to
G(C,C\bb^*,\l_0) = \bb^0$, almost surely.

\noindent {\em (ii)} Corollary 1 implies that under the given
conditions, the Dantzig selector is consistent for $\bb^*$ if and only
if $\l_n \to 0$.  Furthermore, it gives the almost sure limit of
$\hat{\bb}^{ds}$ in cases where the Dantzig selector is not
consistent (that is, when $\l_n \to \l_0 > 0$).  

\vspace{.1in}

\noindent {\bf Corollary 2. }{\em  Suppose that $E(\e_i^2) = \s^2 < \infty$.  Also assume that $n^{-1}X^TX \to C \in \mathscr{P}_0^+$, that $\lim_{n \to \infty}
n^{-1}\max_{1 \leq i \leq
  n} ||\x_i||^2 = 0$, and that $\sqrt{n}\l_n \to \tilde{\l}_0$.  Let
$A^* = \{j;\ \b^*_j \neq 0\}$ and let $\bar{A}^*$ denote the complement of
$A^*$ in $\{1,...,p\}$.  Then
$\sqrt{n}(\hat{\bb}^{ds} - \bb^*) \Dto \u^0$, where $\Dto$ denotes
convergence in distribution, $\u^0$ solves the optimization problem
\begin{equation} \label{cor2a}
\begin{array}{lc}
\mbox{minimize} & ||\u_{\bar{A}^*}||_1 + \mathrm{sign}(\bb^*)_{A^*}^T\u_{A^*} \\
\mbox{subject to} & ||C\u - \bv^0||_{\infty} \leq \tilde{\l}_0,
\end{array}
\end{equation}
and $\bv^0 \sim N(0,\s^2C)$.}

\vspace{.1in}

Corollary 2 is proved in the Appendix.

\noindent {\em Remarks (i)} The second moment condition on $\e_i$ and
the condition $ n^{-1}\max_{1 \leq i \leq n} ||\x_i||^2 \to 0$ ensure that
$n^{-1/2}X^T\ee$ is asymptotically normal.

\noindent {\em (ii)} If $\tilde{\l}_0 = 0$, then
$\hat{\bb}^{ds}$ has the same asymptotic distribution as the ordinary
least squares estimator.  If $\tilde{\l}_0 > 0$, then the limiting
distribution of the Dantzig selector is not normal.

\noindent {\em (iii)} Corollary 2 should be compared
with Theorem 2 of \citep{knight2000asymptotics}, which describes the
limiting distribution of $\hat{\bb}^{lasso}$.  Though the limiting
distribution of lasso is determined by an unconstrained
optimization problem, the term $||\u_{\bar{A}^*}||_1 + \mathrm{sign}(\bb^*)_{A^*}^T
\u_{A^*}$ in the limiting  optimization problem for the Dantzig
selector (\ref{cor2a}) also appears in the limiting optimization problem for lasso.

\section{Discussion}

The results in this paper address fairly long-standing open
questions about uniqueness for the Dantzig selector and lasso.  To
summarize, we prove that the Dantzig selector and lasso estimators are unique in almost all
instances.  Though these results may appear to be
somewhat esoteric, Proposition 2 and its corollaries demonstrate their
potential usefulness.  Indeed, we have shown that once uniqueness is
understood, it is straightforward to obtain the almost sure limit and
limiting distribution of Dantzig selector estimators.  Taking a
broader view, the results presented here may help clear the path for a more
operator theoretic approach to studying the Dantzig selector, lasso,
and other regularized regression procedures.  Such an approach may
offer additional insights into properties of these methods in a
variety of settings. For instance, one could potential obtain a better
understanding of the Dantzig selector in as asymptotic regime where $p
\to \infty$, which is often of particular interest in regularized
regression problems,
by defining the Dantzig selector operator on an appropriate infinite
dimensional space (analogous
to the operator $G$ defined in Section 4 above) and studying its
continuity properties in this more abstract setting.  Future research in this direction is needed.

\section*{Appendix}

{\em Proof of Proposition 1. } The following two lemmas establishes the
Karush-Kuhn-Tucker (KKT) conditions for the Dantzig selector and lasso
optimization problems.  The lemmas appear in various forms in several
references, including \citep{efron2007discussion}, \citep{asif2008masters},
\citep{dicker2010thesis}, and \citep{asif2010lasso}, and proofs are
omitted.

\vspace{.1in}
\noindent {\bf Lemma A1. } {\em The vector $\hat{\bb} = \hat{\bb}^{ds}\in \R^p$ is a solution
  to the Dantzig selector (\ref{DS}) if and only if there is $\hat{\mm} \in \R^p$
  such that
\begin{eqnarray}
\label{lemmaA1a} n^{-1}||X^T(\y - X\hat{\bb})||_{\infty} & \leq & \l \\
\label{lemmaA1b} n^{-1}||X^TX\hat{\mm}||_{\infty} & \leq & 1 \\
\label{lemmaA1c} n^{-1}\hat{\mm}^TX^TX\hat{\bb} & = & ||\hat{\bb}||_1 \\
\label{lemmaA1d} n^{-1}\hat{\mm}^TX^T(y - X\hat{\bb}) & = & \l ||\hat{\mm}||_1.
\end{eqnarray}}

\vspace{.1in}

\noindent {\bf Lemma A2. } {\em The vector $\hat{\bb} = \hat{\bb}^{lasso} \in
  \R^p$ is a solution to the lasso optimization problem (\ref{lasso})
if and only if
\[
\begin{array}{rcll}
n^{-1}\X_j^T(\y - X\hat{\bb}) & = & \l \mathrm{sign}(\hat{\b}_j) \ \ &
\mbox{if } \hat{\b}_j \neq 0 \\
|n^{-1}\X_j^T(\y - X\hat{\bb})| & \leq & \l &  \mbox{if } \hat{\b}_j = 0.
\end{array}
\]}

To prove 1 (a), we assume that [Mult] holds and show that the Dantzig
selector has multiple solutions.  Let $\hat{\bb} = \bb^0$, $\hat{\mm}
= \mm^0$, and $A,B \subseteq \{1,...,p\}$ be as in [Mult] and take $\u \in \R^p\setminus\{0\}$ so that
$\u_{\bar{A}} = 0$ and $n^{-1}X_B^TX_A\u_A =0$, where $\bar{A}$ is the
complement of $A$ in $\{1,...,p\}$.  Then it is clear from Lemma A1
that $\hat{\bb} = \bb^0$ is a solution to the Dantzig selector.
Furthermore, using Lemma A1, it is easy to check that $\hat{\bb} = \bb^0
+ t\u$ is a solution to the Dantzig selector for $t \in \R$
sufficiently small (take $\hat{\mm} = \mm^0$).

Now suppose that $\bb^1,\bb^2 \in \R^p$ are
distinct solutions to the Dantzig selector and let $\mm^1,\mm^2 \in
\R^p$ be vectors such that $\hat{\bb} = \bb^i$ and $\hat{\mm} = \mm^i$, $i
= 1,2$ satisfy (\ref{lemmaA1a})--(\ref{lemmaA1d}).  Without loss of
generality, assume that $\ss = \mbox{sign}(\bb^1) =
\mbox{sign}(\bb^2)$ and $\rr =
\mbox{sign}(\mm^1) = \mbox{sign}(\mm^2)$, where we
define $\mbox{sign}(\bb)_j = \mbox{sign}(\b_j) = 
\b_j/|\b_j|$ or $0$, according to $\b_j \neq 0$ or $\b_j = 0$, for
$\bb \in \R^p$.  Let $A =
\{j; \ \b_j^1 \neq 0\} = \{j; \ \b_j^2 \neq 0\}$, $B = \{j; \
\mu_j^1 \neq 0\} = \{j; \ \mu_j^2 \neq 0\}$.  Then
(\ref{lemmaA1b})--(\ref{lemmaA1c}) imply that
$||n^{-1}X^TX_B\mm^i||_{\infty} \leq 1$ and $n^{-1}X_A^TX_B\mm^i = \ss \in
\{\pm1\}^{|A|}$, $i = 1,2$.
Additionally, (\ref{lemmaA1d}) implies that $n^{-1}X_B^TX_A\bb^1 =
n^{-1}X_B^TX_A\bb^2$.  Hence,
$\dim\left[\mathrm{null}(n^{-1}X_B^TX_A)\right] > 0$.  It follows that
$n^{-1}X^TX$ is parallel to the $\ell^1$-ball.

Finally, to prove Proposition 1 (c), suppose
\[
\bb^1,\bb^2 \in \argmin_{\bb \in \R^p} \frac{1}{2n} ||y - X\bb||^2 + \l ||\bb||_1
\]
are distinct and suppose without loss of generality that
$\ss = \mbox{sign}(\bb^1) = \mbox{sign}(\bb^2)$.  Let $A
= \{j; \ \bb^1_j \neq 0\} = \{j; \ \bb^2_j \neq 0\}$ and $\u = \bb^2 -
\bb^1$.  Notice that for $0 \leq t \leq 1$ we have
\begin{eqnarray}\nonumber
\frac{1}{2n}\left\{2t\u^TX^T(\y - X\bb^1) -
 t^2\u^TX^TX\u\right\} & = &\frac{1}{2n}\left\{||\y - X\bb^1||^2
\right. \\ \nonumber
&& \qquad \left. - ||\y -
 X(\bb^1 + t\u)||^2\right\} \\ \label{prop1a} & = &
\l t \ss^T\u.
\end{eqnarray}
Since (\ref{prop1a}) must hold for all $0
\leq t \leq 1$ and since $\l > 0$, we must have $X\u = X_A\u_A = 0$ and
$\ss^T\u = 0$.  It follows that
\begin{equation}\label{prop1b}
\mbox{dim}\left[\mbox{null}\left(n^{-1}X^TX_A\right)\right] > 0
\end{equation}
and $t = ||\bb^1||_1 =||\bb^2||_1$.

Now, let $\w = \l^{-1}[(X^TX)^-X^T\y - \bb^1] \in \R^p$, where $(X^TX)^-$ is the
Moore-Penrose pseudoinverse of $X^TX$.  Then Lemma A2 implies
that
\[
n^{-1}||X^TX\w||_{\infty} = \frac{1}{n\l}||X^T(\y - X\bb^1)||_{\infty} \leq 1
\]
and $n^{-1}X_A^TX\w = (n\l)^{-1}X_A^T(\y - X\bb^1) \in \{\pm1\}^{|A|}$.  Proposition 1 (c) follows from these observations plus (\ref{prop1b}).
\hfill $\blacksquare$
\vspace{.1in}

{\em Proof of Proposition 2. } To prove Proposition 2, we make use of the following lemma.

\noindent\textbf{Lemma A3. } {\em 
Suppose that $n \geq p$ and that the rows of $X$ are iid and drawn from a distribution which is continuous with respect to Lebesgue measure on $\R^p$.  Let $W$ be an $n \times q$ matrix of rank $q \leq n$.  Then $X^TW$ has rank $\min\{q,p\}$ with probability 1.}

{\em Proof of Lemma A3. }
Let $X$ and $W$ be as in the statement of the lemma.  Without loss of generality, suppose that $q = p$. When $p = 1$, the result is true.  For $p > 1$, let $[p] = \{1,...,p\}$.  To facilitate a proof by induction, assume that $X_{[p-1]}^TW_{[p-1]}$ has rank $p - 1$ with probability 1.  On the event that $X_{[p-1]}^TW_{[p-1]}$ has rank $p- 1$, the rank of $X^TW$ is less than $p$ if and only if
\begin{equation}\label{p1l}
\X_p^T\left\{I - W_{[p - 1]}(X_{[p-1]}^TW_{[p-1]})^{-1}X_{[p-1]}^T\right\}\W_p = 0,
\end{equation}
where $\X_p = X_{\{p\}}$ and $\W_p = W_{\{p\}}$.  
Since $W$ has full rank, it follows that 
\[
\left\{I - W_{[p -
    1]}(X_{[p-1]}^TW_{[p-1]})^{-1}X_{[p-1]}^T\right\}\W_p \neq
0,
\] with probability 1.  Thus, conditioning on $X_{[p-1]}$ and using
the fact that the conditional distribution of $\X_{p}$ is
continuous, it follows that (\ref{p1l}) holds with probability 0.  We
conclude that $X^TW$ has rank $p$ with probability 1. \hfill $\Box$

Getting back to the proof of Proposition 2, suppose that the rows of $X$ are iid and drawn from a distribution which is continuous with respect to Lebesgue measure on $\R^p$.  Then $X$ has rank $\min\{n,p\}$ with probability 1.  Let $A,B \subseteq \{1,...,p\}$ and decompose $A,B$ so that $A = A_0 \cup J$, $B = B_0 \cup J$, and $A_0$, $B_0$, and $J$ are disjoint.  If $|A| > n$, then $X_B^TX_A$ has a non-trivial null space. Suppose for the moment that $|A| \leq n$.  When $X$ has full rank, the dimension of the null space of $X_B^TX_A$ is non-zero if and only if
\[
\dim \left( \mathrm{null}\left[X_{B_0}^T\left\{I - X_J(X_J^TX_J)^{-1}X_J^T\right\}X_{A_0}\right]\right) > 0.
\]
Furthermore, if $X$ has full rank, then $\left\{I -
  X_J(X_J^TX_J)^{-1}X_J^T\right\}X_{A_0}$ has full rank. Conditioning
on $X_A$ and appealing to Lemma A3, it follows that the rank of $X_{B_0}^T\left\{I - X_J(X_J^TX_J)^{-1}X_J^T\right\}X_{A_0}$ is $\min\{|A_0|,|B_0|\}$ with probability 1.  Thus the null-space of $X_B^TX_A$ is non-trivial with positive probability if and only if $\min\{|B|,n\} < |A|$.

Now suppose that $\min\{|B|,n\} < |A|$.  There are two cases: $|B| < |A| \leq n$ and $n < |A|$.  In each case, the probability that there exists $\w \in \R^{|B|}$ such that $X_A^TX_B\w \in\{\pm1\}^{|A|}$ is 0.  We prove this for the case $|B| < |A| \leq n$; the case $n < |A|$ follows similarly.  Assume that $|B| < |A| \leq n$.  Choose $A_1 \subseteq A_0$ such that $|A_1| = |B_0|$ and let $\tilde{A} = J \cup A_1$.  Suppose that $X_{\tilde{A}}^TX_B\w = \mathbf{s}$ for some $\mathbf{s} \in \{\pm1\}^{|\tilde{A}|}$ and $\w \in \R^{|B|}$.  Then, assuming that $X$ is full rank,
\[
\w_J = (X_J^TX_J)^{-1}(\mathbf{s}_J - X_J^TX_{B_0}\w_{B_0})
\]
and
\[
X_{A_1}^T\left[\left\{I - X_J(X_J^TX_J)^{-1}X_J^T\right\}X_{B_0}\w_{B_0} + X_J(X_J^TX_J)^{-1}\mathbf{s}_J\right] = \mathbf{s}_{A_1}.
\]
Thus, we have
\[
\w_{B_0} = \left[X_{A_1}^T\left\{I - X_J(X_J^TX_J)^{-1}X_J^T\right\}X_{B_0}\right]^{-1}\left\{\mathbf{s}_{A_1} - X_{A_1}^TX_J(X_J^TX_J)^{-1}\mathbf{s}_J\right\},
\]
where Lemma A3 guarantees that $X_{A_1}^T\left\{I -
  X_J(X_J^TX_J)^{-1}X_J^T\right\}X_{B_0}$ is invertible with
probability 1.  Since, conditional on $X_{A_1 \cup B}$, the rows of $X_{A_0\setminus
  A_1}$ are independent and have continuous
distributions with respect to Lebesgue
measure on $\R^{|A_0\setminus A_1|}$, it follows that 
\[
X_{A_0\setminus A_1}^T\!\left[\left\{I \!-\!
    X_J(X_J^TX_J)^{-1}X_J^T\right\}X_{B_0}\w_{B_0} \!+\!
  X_J(X_J^TX_J)^{-1}\mathbf{s}_J\right] \! \in\! \{\pm1\}^{|A_0| - |A_1|}
\]
with probability 0.  Thus, as claimed, the probability that there exists $\w \in \R^{|B|}$ such that $X_A^TX_B\w \in\{\pm1\}^{|A|}$ is 0.

The results from the last two paragraphs imply that
\[
P\left[\begin{array}{c} \dim\left\{ \mathrm{null}(X_B^TX_A)\right\} > 0 \mbox{ and } \\
X_A^TX_B\w \in \{\pm1\}^{|A|} \mbox{ for some} \ \w \in \R^{|B|} \end{array} \right] = 0
\]
It follows that $X^TX$ is parallel to the $\ell^1$-ball with probability 0, as was to be shown.
\hfill $\blacksquare$

{\em Proof of Proposition 3. } For $n \in \mathbb{N}$, let $C_n,C \in
\mathscr{P}_0^+$, $\bv_n,\bv \in \R^p$, and $\l_n,\l \geq 0$ and assume
that $C_n \to C$, $\bv_n \to \bv$, and $\l_n \to \l$.  Let $\hat{\u}_n =
G(C_n,\bv_n,\l_n)$ and let $\hat{\u} = G(C,\bv,\l)$.  We show that
$\hat{\u}_n \to \hat{\u}$.

Since $\sup_n ||\hat{\u}_n|| < \infty$, there exists a subsequence
$\{\hat{\u}_{n_k}\}_{k = 1}^{\infty}$ and a vector $\u_0 \in \R^p$ such
that $\hat{\u}_{n_k} \to \u_0$.  To prove the proposition, it suffices
to show that $\u_0 = \hat{\u}$.  By continuity of the
$\ell^{\infty}$-norm, we must have
\[
||C\u_0 - \bv||_{\infty} \leq \l.
\]
Also, by the optimality properties of $\hat{\u}_{n_k}$, we must have
\begin{equation} \label{prop3a}
||\u_0||_1 = \lim_{k \to \infty} ||\hat{\u}_{n_k}||_1 \leq \liminf_{k \to \infty} ||\w_{n_k}||_1
\end{equation}
for any sequence $\{\w_{n_k}\}$, with $\w_{n_k} \in \R^p$ and
\begin{equation} \label{prop3b}
||C_{n_k}\w_{n_k} - \bv_{n_k}||_{\infty} \leq \l_{n_k}.
\end{equation}
We consider two cases: $\l = 0$ and $\l > 0$.  First suppose $\l = 0$ and define
$\w_{n_k} = C_{n_k}^{-1}(C\hat{\u} + \bv_{n_k} - \bv)$.  Then (\ref{prop3b})
holds and $\w_{n_k} \to \hat{\u}$.  From (\ref{prop3a}), it follows that
$||\u_0||_1 \leq ||\hat{\u}||_1$ and the optimality of $\hat{\u}$ implies
that $\hat{\u} = \u_0$.  Now suppose that $\l > 0$ and define $\w_{n_k} =
(\l_{n_k}/\l)C_{n_k}^{-1}(C\hat{\u} -\bv) +  C_{n_k}^{-1}\bv_{n_k} $.  Then
(\ref{prop3b}) holds and, as in the previous case, we conclude that
$\hat{\u} = \u_0$.   Thus, in either case, $\hat{\u} = \u_0$, as was to be
shown. \hfill $\blacksquare$

\vspace{.1in}

{\em Proof of Corollary 2. }  The conditions $E(\e_i^2) <
\infty$ and $n^{-1} \max_{1 \leq i \leq n} ||\x_i||^2 \to 0$ ensure
that $n^{-1/2}X^T\ee \Dto \bv^0 \sim N(0,\s^2C)$, by the Lindeberg-Feller central
limit theorem.  By the Skorokhod
representation theorem, we may assume without loss of generality that
$n^{-1/2}X^T\ee \to \bv^0$ almost surely.

Now let $\u = \sqrt{n}(\bb - \bb^*)$ and notice that the Dantzig selector
(\ref{DS}) is equivalent to the optimization problem
\begin{equation}\label{cor2b}
\begin{array}{lc}
\mbox{minimize} & ||\sqrt{n}\bb^* + \u||_1 \\
\mbox{subject to} & \left|\left|n^{-1}X^TX\u -
    n^{-1/2}X^T\ee\right|\right|_{\infty} \leq \sqrt{n} \l_n.
\end{array}
\end{equation}
In particular, $\hat{\u} = \sqrt{n}(\hat{\bb}^{ds} - \bb^*)$ solves
(\ref{cor2b}).  We show that $\hat{\u} \to \u^0$, the solution to
(\ref{cor2a}), almost surely.  This suffices to prove the corollary.

Since $n^{-1}X^TX \to C$, $n^{-1/2}X^T\ee \to \bv^0$ almost surely, and
$\sqrt{n}\l_n \to \tilde{\l}_0$, it
follows that there is an almost surely finite random variable $M$ such
that $||\u||_{\infty} \leq M/2$ whenever $\u$ is feasible for the optimization
problem (\ref{cor2b}).  Let $\ss = \mathrm{sign}(\bb^*)$ and notice
that if $\sqrt{n}\min \{|\b_j^*|; \ j \in A^*\} > M$ and $\u$ is
feasible for (\ref{cor2b}), then $\mathrm{sign}(\sqrt{n}\bb^* + \u) = \mathrm{sign}(M\ss +
\u)$.  It follows that
\[
G(n^{-1}X^TX,n^{-1/2}X^T\ee + n^{-1}X^TXM\ss,\sqrt{n}\l_n) = M\ss +
\hat{\u}
\]
whenever $\sqrt{n}\min \{|\b_j^*|; \ j \in A^*\} > M$.  Taking $n \to
\infty$, Proposition 3 implies that $\hat{\u} \to G(C,\bv^0 +
CM\ss,\tilde{\l}_0) - M\ss$ almost surely and it is straightforward to
check that $\u^0 = G(C,\bv^0 +
CM\ss,\tilde{\l}_0) - M\ss$. \hfill $\blacksquare$

\end{document}